А. Е. Максимов, С. П. Куксенко

# Эффективный расчет каузальной емкостной матрицы многопроводной линии передачи в диапазоне изменения ее параметров методом моментов

**Аннотация**

В работе предложен метод неравномерной сегментации микрополосковых многопроводных линий передачи с торцевой связью, который позволяет получить физичные емкостные матрицы с наименьшими вычислительными затратами при многовариантном анализе в диапазоне изменения параметров линий. Показано, что использование этого метода позволяет достичь экономии времени до 49%.

**Введение**

Конкуренция между производителями радиоэлектронных средств (РЭС) вызывает необходимость регулярного обновления и совершенствования всех их видов. При этом с ростом сложности РЭС, их разработка и совершенствование становятся невозможными без применения автоматизированного проектирования, в основе которого лежит компьютерное моделирование [1]. Для минимизации вычислительных затрат на разработку РЭС необходимо уменьшать вычислительную сложность методов компьютерного моделирования при контролируемой точности результатов.

Известно, что основными элементами РЭС являются линии передачи (ЛП) [2], для моделирования которых широко применяется квазистатический подход [3], который основан на решении уравнения Пуассона/Лапласа с помощью численных методов конечных разностей (finite difference method, FDM) [4], конечных элементов (finite element method, FEM) [5] или моментов (method of moments, MoM) [6]. При этом в MoM используется поверхностная сегментация анализируемой ЛП, что априорно снижает вычислительную сложность решения. При использовании MoM решается уравнение Пуассона/Лапласа в интегральном виде [7]

$$\varphi(\mathbf{r}) = \frac{1}{\varepsilon_0} \int \sigma(\mathbf{r}')G(\mathbf{r},\mathbf{r}')d\Gamma, \qquad (1)$$

где $\sigma(\mathbf{r}')$ – поверхностная плотность заряда, $\mathbf{r}$ и $\mathbf{r}'$ – точки наблюдения $(x, y)$ и расположения источника $(x', y')$, соответственно, $d\Gamma$ – дифференциал по поверхности ЛП, $G(\mathbf{r}, \mathbf{r}')$ – функция Грина, а $\varepsilon_0$ – электрическая постоянная. Граничные условия задаются по приложенному напряжению φ и требуется найти плотность заряда σ. Это реализуется за счет сведения уравнения (1) к системе линейных алгебраических уравнений (СЛАУ), решаемой как прямыми [8, 9], так и итерационными [10, 11] методами. Затем из решения СЛАУ вычисляются матрицы погонных параметров ЛП **C**, **G**, **R** и **L** [12 – 14], вычисление первой из которых (емкостной матрицы) представляет наибольшую сложность [15].

При этом на практике, при анализе некоторых многопроводных линий передачи (МПЛП) матрицы погонных параметров могут оказаться нефизичными [15 – 17], что приводит к недопустимой нефизичности отклика в конце ЛП [18, 19]. Наиболее простым и разумным методом борьбы с нефизичностью является учащение поверхностной сегментации границ ЛП, однако часто это приводит к росту ошибок округления [4, 20].

Методы сегментации можно условно разбить на две группы: равномерная, когда все границы ЛП разбиваются на сегменты одинаковой длины, и неравномерная, когда сегменты имеют разную длину, которая может задаваться как вручную, так и с применением различных адаптивных алгоритмов [21, 22]. Общий принцип адаптивного учащения состоит в том, что учащение сегментации в областях сосредоточения особенностей решения (например, торцы проводников при торцевой связи) в наибольшей степени повышает точность результатов [23]. В данной работе для поиска метода устранения нефизичности результатов исследована как равномерная, так и неравномерная адаптивная сегментация.

Кроме того, выполнен поиск эффективного средства уменьшения вычислительных затрат на многовариантный анализ МПЛП, часто используемый на практике на начальных этапах проектирования РЭС [24, 25]. Для этого выполнен учет частичных изменений в матрице СЛАУ. Так, при анализе МПЛП в диапазоне изменения ее параметров свое положение в пространстве меняют лишь некоторые границы ЛП, поэтому нет необходимости пересчитывать матрицу целиком.

Цель работы – поиск оптимального метода сегментации МПЛП, дающего физичные результаты с наименьшими вычислительными затратами при многовариантном анализе МПЛП.

Работа построена следующим образом. В разделе 1 кратко рассмотрены основные положения математической модели для вычисления емкостной матрицы **C**, в разделе 2 – применение равномерной сегментации. Раздел 3 посвящен разработке метода неравномерной сегментации для устранения нефизичного вычисления емкостной матрицы **C**, а раздел 4 – сокращению вычислительных затрат на многовариантный анализ МПЛП с использованием этого метода.

## 1. Математическая модель вычисления матрицы C

Математические модели для вычисления матрицы **C** с использованием MoM для двумерных структур с границами произвольной сложности хорошо разработаны в [14]. Для ясности дальнейшего изложения кратко поясним суть этих моделей.

Когда проводники ЛП соприкасаются с диэлектриками (границы проводник-диэлектрик), вычисления выполняются с учетом полной плотности заряда $\sigma_T$, которая представляет собой сумму плотности свободного $\sigma_S$ и поляризационного $\sigma_P$ зарядов: $\sigma_T(\mathbf{r}) = \sigma_S(\mathbf{r}) + \sigma_P(\mathbf{r})$ [15]. При этом на границе диэлектрик-диэлектрик полная плотность заряда равна $\sigma_P$. Тогда, для границ проводник-диэлектрик, (1) преобразуется к виду

$$\varphi(\mathbf{r}) = -\frac{1}{2\pi\varepsilon_0} \int_{L_C} \sigma_T(\mathbf{r}') \ln|\mathbf{r} - \mathbf{r}'| \, dl', \; \mathbf{r} \in L_C, \quad (2)$$

где $dl'$ – элемент контура границ, $L_C$ – контур проводниковых границ, а для границ диэлектрик-диэлектрик –

$$0 = \frac{1}{2\pi\varepsilon_0} \int_{L_D} \sigma_T(\mathbf{r}') \frac{\mathbf{r} - \mathbf{r}'}{|\mathbf{r} - \mathbf{r}'|^2} \cdot \mathbf{n} \, dl' - \frac{\varepsilon_{r2} + \varepsilon_{r1}}{\varepsilon_{r2} - \varepsilon_{r1}} \frac{\sigma_T(\mathbf{r})}{2\varepsilon_0}, \; \mathbf{r} \in L_D, \quad (3)$$

где $L_D$ – контур диэлектрических границ, $\mathbf{n}$ – вектор нормали, а $\varepsilon_{r1}$ и $\varepsilon_{r2}$ – относительные диэлектрические проницаемости среды на положительной (на которую направлен $\mathbf{n}$) и отрицательной (от которой направлен $\mathbf{n}$) сторонах границы соответственно [14].

Используя поверхностную сегментацию границ МПЛП, а также базисные и тестовые функции, уравнения (2) и (3) преобразуются к матричному уравнению вида $\mathbf{S\Sigma} = \mathbf{V}$, где $\mathbf{S}$ – плотная матрица порядка $N$, $\mathbf{\Sigma}$ и $\mathbf{V}$ – матрицы размера $N \times N_C$, $N_C$ – число проводников в МПЛП, не считая опорного. Элементы $\mathbf{S}$ вычисляются аналитически [14], а решение СЛАУ $N_C$ раз, с каждым столбцом $\mathbf{v}$ из $\mathbf{V}$ в качестве вектора свободных членов, даёт распределение заряда на поверхности МПЛП, из которого затем вычисляется емкостная матрица **C**, выражающая связь между зарядами и напряжениями на проводниках: $\mathbf{q} = \mathbf{Cv}$, где $\mathbf{q}$ и $\mathbf{v}$ – вектор-столбцы заряда и напряжения соответственно.

Полученная матрица **C** для МПЛП с $N_C$ проводниками имеет вид [26]

$$\mathbf{C} = \begin{bmatrix} c_{11} + c_{12} + \ldots + c_{1N_C} & -c_{12} & \ldots & -c_{1N_C} \\ -c_{21} & c_{21} + c_{22} + \ldots + c_{2N_C} & \ldots & -c_{2N_C} \\ \ldots & \ldots & \ldots & \ldots \\ -c_{N_C 1} & -c_{N_C 2} & \ldots & c_{N_C 1} + c_{N_C 2} + \ldots + c_{N_C N_C} \end{bmatrix}. \quad (4)$$

Из (4) видно, что матрица **C** должна быть симметричной, а сумма абсолютных значений внедиагональных элементов в каждой ее строке не должна превышать значение диагонального элемента. Однако эти условия не всегда выполняются [18, 19]. Так, усложнение конфигурации МПЛП за счет увеличения числа проводников (сигнальных, опорных, и т. д.) приводит к увеличению порядка матрицы **C** и к её нефизичности. Часто такое наблюдается в микрополосковых МПЛП с числом проводников более восьми, имеющих торцевые связи (рисунок 1).

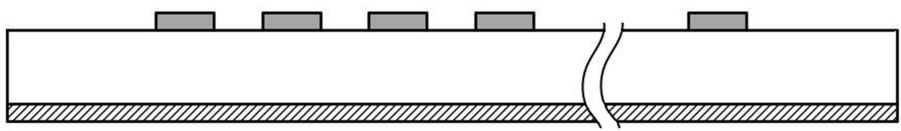

Рисунок 1 – Поперечное сечение типовой микрополосковой МПЛП

Нефизичность матрицы **C** в таких структурах может также проявляться в положительности её внедиагональных элементов, в увеличении взаимной емкости проводников при увеличении расстояния между ними, а также в асимметрии матрицы **C** [15 – 17].

## 2. Равномерная сегментация границ ЛП

Рассмотрим метод равномерной сегментации границ ЛП. В [27] приведены примеры вычисления емкостной матрицы с достаточной точностью при использовании равномерной сегментации с длиной сегмента $t/3$, где $t$ – толщина проводников МПЛП. Однако с ростом $t$ такая сегментация часто дает недостаточно точные результаты.

Продемонстрируем это на примере 8-проводной ЛП с торцевой связью, поперечное сечение которой приведено на рисунке 2. Параметры ЛП: ширина проводников $w = 0{,}05$ мм; расстояние между проводниками $s = 0{,}05$ мм; расстояние от крайнего проводника до внешней границы ЛП $d = 0{,}15$ мм; толщины диэлектриков $h_1 = h_3 = 0{,}05$ мм; $h_2 = 0{,}15$ мм; относительные диэлектрические проницаемости $\varepsilon_{r1} = \varepsilon_{r3} = 3{,}8$, $\varepsilon_{r2} = 2$.

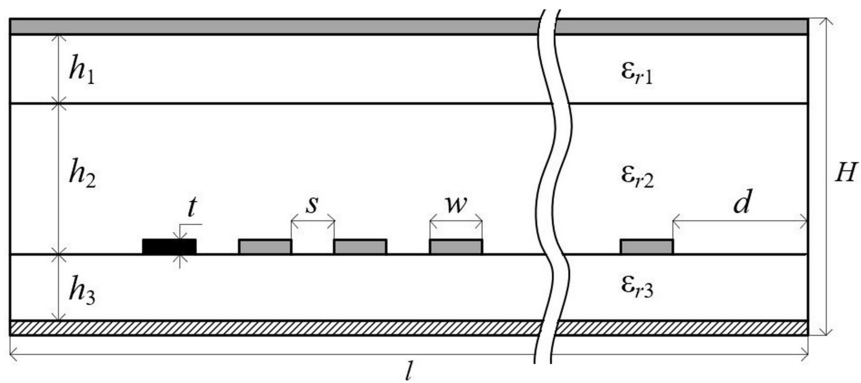

Рисунок 2 – МПЛП 1: поперечное сечение

Использовано 5 типовых значений $t$: 0,005; 0,018; 0,05; 0,07 и 0,105 мм. Выполнена оценка результатов вычисления матрицы **C** при сегментации $t/n$, где $n = 1, 3, 5, 7$ и 9.

В качестве контролируемых использованы элементы первой строки матрицы **C**, т. е. коэффициенты $C_{11}$ контрольного проводника (черный проводник на рисунке 2) и $C_{12}$ – $C_{18}$. Толщина проводников $t$ изменялась таким образом, чтобы общие размеры МПЛП ($l$ и $H$) не менялись. При вычислениях использована реализация математических моделей, описанных в разделе 1 в пакете GNU Octave.

В таблицу 1 сведены полученные значения контролируемых величин при $t = 0{,}005$ мм. Аналогичные результаты для остальных $t$ приведены в Приложении А.

Таблица 1 – Значения (пФ/м) первой строки матрицы **C** при $t = 0{,}005$ мм и изменении $n$, полученные в GNU Octave для МПЛП 1 (рисунок 2)

| $n$ | $C_{11}$ | $|C_{12}|$ | $|C_{13}|$ | $|C_{14}|$ | $|C_{15}|$ | $|C_{16}|$ | $|C_{17}|$ | $|C_{18}|$ |
|---|---|---|---|---|---|---|---|---|
| 1 | 89,00 | 12,80 | 0,64 | 9,89e–2 | 1,90e–2 | 4,80e–3 | 1,91e–3 | 1,25e–3 |
| 3 | 92,10 | 12,40 | 0,63 | 9,46e–2 | 1,62e–2 | 2,82e–3 | 5,02e–4 | 1,44e–4 |
| 5 | 92,00 | 12,50 | 0,63 | 9,52e–2 | 1,63e–2 | 2,84e–3 | 5,05e–4 | 1,44e–4 |
| 7 | 92,00 | 12,50 | 0,63 | 9,54e–2 | 1,63e–2 | 2,85e–3 | 5,07e–4 | 1,44e–4 |
| 9 | 92,00 | 12,50 | 0,63 | 9,55e–2 | 1,63e–2 | 2,85e–3 | 5,07e–4 | 1,44e–4 |

Различия в значениях матрицы **C** оценивались по контролируемой величине

$$\Delta_C = \frac{\left|C_{t/(n-2)} - C_{t/n}\right|}{C_{t/n}}, \qquad (5)$$

где $n = 3, 5, 7, 9$. В таблице 2 приведены различия по $C_{11}$. Различия в матричных нормах, оцениваемые как

$$\Delta_F = \frac{\left\|\mathbf{C}_{t/(n-2)} - \mathbf{C}_{t/n}\right\|_F}{\left\|\mathbf{C}_{t/n}\right\|_F}, \qquad (6)$$

где $\|\cdot\|_F$ – норма Фробениуса [10], сведены в таблицу 3.

Таблица 2 – Различия по $\Delta_{C_{11}}$ для МПЛП 1 (рисунок 2), %

| $n$ | $t$, мм | | | | |
|---|---|---|---|---|---|
| | 0,005 | 0,018 | 0,05 | 0,07 | 0,105 |
| 3 | 3,32 | 0,6 | 0,7 | 1,46 | 2,89 |
| 5 | 0,03 | 0,08 | 0,03 | 1,16 | 2,69 |
| 7 | 0,02 | 0,06 | 0,05 | 0,31 | 0,44 |
| 9 | 0,01 | 0,01 | 0,04 | 0,15 | 0,24 |

Таблица 3 – Различия по $\Delta_F$ для МПЛП 1 (рисунок 2), %

| $n$ | $t$, мм | | | | |
|---|---|---|---|---|---|
| | 0,005 | 0,018 | 0,05 | 0,07 | 0,105 |
| 3 | 3,08 | 1,35 | 2,78 | 2,35 | 12,5 |
| 5 | 0,06 | 0,22 | 0,43 | 1,37 | 2,13 |
| 7 | 0,03 | 0,07 | 0,14 | 0,31 | 0,44 |
| 9 | 0,01 | 0,03 | 0,07 | 0,15 | 0,19 |

Из таблиц 2 и 3 видно, что при $t \leq 0{,}05$ мм переход от сегментации $t/3$ к сегментации $t/5$ дает различие менее 1%, и поэтому в этом случае целесообразно использовать сегментацию $t/3$, а при $t > 0{,}05$ мм – $t/5$, т. к. переход от $t/3$ к $t/5$ дает различие более 1%, а от $t/5$ к $t/7$ – менее 1%. Стоит отметить, что в таблице 2 для $t = 0{,}018$ мм и $t = 0{,}05$ мм различие менее 1% наблюдается и при переходе от сегментации $t/1$ к сегментации $t/3$, однако, сегментацию менее чем $t/3$ брать не целесообразно, потому что она плохо аппроксимирует заряд на торцах проводников, а значит, результаты априорно будут неточны. В исследовании рассматривается различие менее 1% аналогично работе [28].

### 3. Некаузальное поведение и неравномерная сегментация границ ЛП

Даже при частой равномерной сегментации ($t/9$), как было показано в предыдущем разделе, в некоторых случаях (таблицы А.2 – А.4) наблюдается нефизичное поведение, проявляющееся в увеличении расчетной емкостной связи с увеличением расстояния между проводниками. Более того, при увеличении числа проводников растет нефизичность результатов. Так, из таблицы 4 видно, что абсолютные значения $C_{19}$ и $C_{110}$ превышают значение $C_{18}$, чего априорно не должно быть для данного типа МПЛП.

Таблица 4 – Значения (пФ/м) первой строки матрицы **C** при $t = 0{,}005$ мм и сегментации $t/3$, полученные в GNU Octave для МПЛП 1 (рисунок 2)

| $C_{11}$ | $C_{12}$ | $C_{13}$ | $C_{14}$ | $C_{15}$ | $C_{16}$ | $C_{17}$ | $C_{18}$ | $C_{19}$ | $C_{110}$ |
|---|---|---|---|---|---|---|---|---|---|
| 98,20 | –10,00 | –0,14 | –5,06e–3 | –2,09e–4 | –8,80e–6 | –4,21e–7 | –6,40e–8 | –8,74e–8 | –1,98e–6 |

Выполнено сравнение полученных результатов с другими программами. Использованы TALGAT [29] и FasterCap [30]. В TALGAT использована равномерная сегментация $t/3$, а в FasterCap – адаптивная сегментация по алгоритму учащения Аппеля [31]. Полученные результаты сведены в таблицы 5 и 6 соответственно.

Таблица 5 – Значения (пФ/м) первой строки матрицы **C** при $t = 0{,}005$ мм и сегментации $t/3$, полученные в TALGAT для МПЛП 1 (рисунок 2)

| $C_{11}$ | $C_{12}$ | $C_{13}$ | $C_{14}$ | $C_{15}$ | $C_{16}$ | $C_{17}$ | $C_{18}$ | $C_{19}$ | $C_{110}$ |
|---|---|---|---|---|---|---|---|---|---|
| 97,70 | –10,00 | –0,14 | –5,15e–3 | –2,14e–4 | –8,98e–6 | –3,78e–7 | –1,76e–8 | –4,27e–8 | –1,78e–6 |

Таблица 6 – Значения (пФ/м) первой строки матрицы **C** при $t = 0{,}005$ мм и адаптивной сегментации, полученные в FasterCap (стандартные настройки) для МПЛП 1 (рисунок 2)

| $C_{11}$ | $C_{12}$ | $C_{13}$ | $C_{14}$ | $C_{15}$ | $C_{16}$ | $C_{17}$ | $C_{18}$ | $C_{19}$ | $C_{110}$ |
|---|---|---|---|---|---|---|---|---|---|
| 96,70 | –9,70 | –0,30 | 3,23e–2 | 1,44e–2 | 4,87e–3 | 5,29e–3 | –2,03e–3 | –2,03e–2 | 1,67e–2 |

Из таблиц 5 и 6 видно, что в TALGAT нефизичное поведение сохраняется для тех же элементов, а в FasterCap – для $C_{17} - C_{110}$. Кроме того, в FasterCap элементы $C_{14} - C_{17}$ и $C_{110}$ имеют положительные значения, что усиливает нефизичность результатов. Стоит отметить, что результаты для FasterCap получены при стандартных настройках программы. Однако изменение различных конфигурационных параметров (relative error, mesh relative refinement value, GMRES iteration tolerance, mesh refinement ratio, activating the use Galerkin scheme) в широких пределах не позволило получить физичных результатов.

Далее исследовано использование неравномерной адаптивной сегментации. В работе [22] предложен метод для анализа ЛП методом моментов. Однако этот метод слабо применим для анализа МПЛП ввиду необходимости решения СЛАУ для каждого коэффициента треугольной емкостной матрицы. В работе [27] этот метод был модифицирован на случай МПЛП (множество СЛАУ заменено матричным уравнением). Поясним его суть с помощью упрощенного псевдокода:

1   Задать исходные параметры МПЛП, требуемую точность вычислений *tol*, начальную сегментацию, максимальное число итераций для решения матричного уравнения ($N_{it}^{max}$) и определить контролируемую величину $K$.
2   Вычислить значение $K$.
3   *Для i от 1 до $N_{it}^{max}$*
4       Участить 25% от общего числа сегментов с максимальными значениями плотности заряда на них.
5       Вычислить значение $K$.
6       *Если* $tol < \|K_i - K_{i-1}\|/\|K_{i-1}\|$
7           Продолжать итерации.
8       *Иначе*
9           Выйти из цикла.

Однако и этот метод не решает проблему. В таблицу 7 сведены результаты использования этого метода для 10-проводной ЛП. Видно, что он также дает нефизичные результаты.

Таблица 7 – Значения (пФ/м) первой строки матрицы **C** при $t = 0,005$ мм, полученные в GNU Octave для МПЛП 1 (рисунок 2)

| $C_{11}$ | $C_{12}$ | $C_{13}$ | $C_{14}$ | $C_{15}$ | $C_{16}$ | $C_{17}$ | $C_{18}$ | $C_{19}$ | $C_{110}$ |
|---|---|---|---|---|---|---|---|---|---|
| 97,60 | –9,90 | –0,14 | –4,87e–3 | –2,00e–4 | –9,16e–6 | –1,07e–6 | –6,05e–7 | –4,48e–7 | –1,98e–6 |

Поскольку нефизичным оказался лишь элемент $C_{110}$, то выполнена модернизация метода таким образом, чтобы учащалось не 25% от общего числа сегментов, а $k$% для каждой правой части матричного уравнения (далее метод I). В результате, разбиению подлежат не только сегменты, имеющие наибольший заряд в решении для первого проводника, а учитываются решения для всех проводников.

Поскольку за одну итерацию не может быть разбито более 100% границ ЛП, то процент разбиваемых сегментов не должен превышать значения 100/$N$. В предварительных вычислительных экспериментах установлено, что при $k$ в диапазоне 65 – 85 результирующая матрица **C** физична. Поэтому далее принято $k = 75$. Полученная по методу I матрица **C** (пФ/м) для 10-проводной ЛП имеет вид

$$\begin{bmatrix}
98{,}48 & -9{,}95 & -0{,}14 & -5{,}11\mathrm{e}-3 & -3{,}22\mathrm{e}-4 & -8{,}37\mathrm{e}-5 & -5{,}00\mathrm{e}-5 & -3{,}45\mathrm{e}-5 & -2{,}49\mathrm{e}-5 & -2{,}44\mathrm{e}-5 \\
-9{,}95 & 100{,}23 & -9{,}93 & -0{,}14 & -5{,}10\mathrm{e}-3 & -3{,}28\mathrm{e}-4 & -8{,}72\mathrm{e}-5 & -5{,}14\mathrm{e}-5 & -3{,}46\mathrm{e}-5 & -3{,}32\mathrm{e}-5 \\
-0{,}14 & -9{,}93 & 100{,}23 & -9{,}93 & -0{,}14 & -5{,}10\mathrm{e}-3 & -3{,}28\mathrm{e}-4 & -8{,}72\mathrm{e}-5 & -5{,}16\mathrm{e}-5 & -4{,}29\mathrm{e}-5 \\
-5{,}12\mathrm{e}-3 & -0{,}14 & -9{,}93 & 100{,}23 & -9{,}93 & -0{,}14 & -5{,}10\mathrm{e}-3 & -3{,}28\mathrm{e}-4 & -8{,}74\mathrm{e}-5 & -5{,}86\mathrm{e}-5 \\
-3{,}32\mathrm{e}-4 & -5{,}10\mathrm{e}-3 & -0{,}14 & -9{,}93 & 100{,}23 & -9{,}93 & -0{,}14 & -5{,}10\mathrm{e}-3 & -3{,}28\mathrm{e}-4 & -9{,}26\mathrm{e}-5 \\
-9{,}26\mathrm{e}-5 & -3{,}28\mathrm{e}-4 & -5{,}10\mathrm{e}-3 & -0{,}14 & -9{,}93 & 100{,}23 & -9{,}93 & -0{,}14 & -5{,}10\mathrm{e}-3 & -3{,}32\mathrm{e}-4 \\
-5{,}86\mathrm{e}-5 & -8{,}74\mathrm{e}-5 & -3{,}28\mathrm{e}-4 & -5{,}10\mathrm{e}-3 & -0{,}14 & -9{,}93 & 100{,}23 & -9{,}93 & -0{,}14 & -5{,}12\mathrm{e}-3 \\
-4{,}29\mathrm{e}-5 & -5{,}16\mathrm{e}-5 & -8{,}72\mathrm{e}-5 & -3{,}28\mathrm{e}-4 & -5{,}10\mathrm{e}-3 & -0{,}14 & -9{,}93 & 100{,}23 & -9{,}93 & -0{,}14 \\
-3{,}32\mathrm{e}-5 & -3{,}46\mathrm{e}-5 & -5{,}14\mathrm{e}-5 & -8{,}72\mathrm{e}-5 & -3{,}28\mathrm{e}-4 & -5{,}10\mathrm{e}-3 & -0{,}14 & -9{,}93 & 100{,}23 & -9{,}95 \\
-2{,}44\mathrm{e}-5 & -2{,}49\mathrm{e}-5 & -3{,}45\mathrm{e}-5 & -5{,}00\mathrm{e}-5 & -8{,}37\mathrm{e}-5 & -3{,}22\mathrm{e}-4 & -5{,}11\mathrm{e}-3 & -0{,}14 & -9{,}95 & 98{,}48
\end{bmatrix}.$$

Дополнительная верификация метода I произведена на 10-проводной ЛП из [32], поперечное сечение которой приведено на рисунке 3. Геометрические параметры ЛП (мм): $t = 0{,}02$; $w_1 = w_9 = 0{,}2$; $w_2 = w_8 = w_{10} = 0{,}3$; $w_3 = w_7 = 0{,}4$; $w_4 = w_6 = 0{,}5$; $w_5 = 0{,}6$; $s_1 = s_4 = s_6 = s_9 = 0{,}25$; $s_2 = s_7 = 0{,}3$; $s_3 = s_8 = 0{,}35$; $s_5 = 0{,}2$; $d = 2{,}48$; $h = 1$. Относительная диэлектрическая проницаемость $\varepsilon_r = 4$. В таблицу 8 сведены полученные значения контрольных величин. Из таблицы видно, что все элементы физичны.

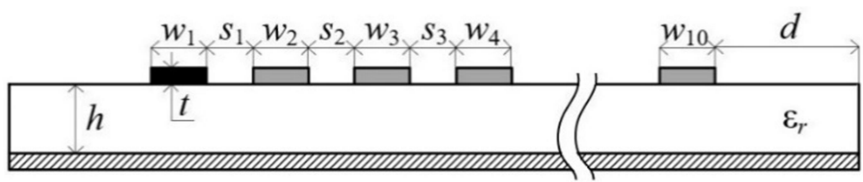

Рисунок 3 – МПЛП 2: поперечное сечение

Таблица 8 – Значения (пФ/м) первой строки матрицы **C**, полученные в GNU Octave для МПЛП 2 (рисунок 3)

| $C_{11}$ | $C_{12}$ | $C_{13}$ | $C_{14}$ | $C_{15}$ | $C_{16}$ | $C_{17}$ | $C_{18}$ | $C_{19}$ | $C_{110}$ |
|---|---|---|---|---|---|---|---|---|---|
| 49,50 | –19,00 | –2,40 | –0,58 | –0,26 | –0,15 | –0,11 | –9,97e–2 | –8,69e–2 | –7,09e–2 |

### 4. Многовариантный анализ

При переходе от одновариантного анализа к многовариантному возникает необходимость многократно переформировывать матрицу СЛАУ. При этом в матрице изменяются лишь те элементы, которые соответствуют сместившимся в пространстве сегментам МПЛП. Кроме того, при явном задании числа сегментов для каждой физической границы ЛП, их общее число остается неизменным. Поэтому, при изменении геометрических параметров МПЛП в широком диапазоне значений, соответствующие элементы матрицы СЛАУ будут располагаться на тех же позициях. Именно этот факт позволяет при многовариантном анализе учитывать лишь частичные изменения в матрице СЛАУ.

Стоит отметить, что для достижения наибольшей экономии вычислительных ресурсов важно выполнять геометрическое построение ЛП таким образом, чтобы при изменении выбранного параметра как можно меньшее число сегментов меняло свое положение в пространстве. Это правило не требуется выполнять, если изменяемый параметр не относится к геометрическим (например, при изменении относительной диэлектрической проницаемости). Также целесообразно зафиксировать геометрические размеры МПЛП, и, например, изменение $w$ и $s$ компенсировать изменением $d$ (как и сделано в работе). Далее приведен псевдокод метода I для многовариантного анализа МПЛП при изменении ее геометрических параметров:

1 Задать параметры ЛП, в том числе число значений изменяемых параметров $N_{par}$ и процент учащаемых сегментов $k$.

2 Задать начальную сегментацию.

3 Сформировать и решить матричное уравнение (произвести предварительное вычисление).

4 Определить подлежащие разбиению сегменты.

5 *Для i от 1 до $N_{par}$*

6    *Если i = 1*

7       Выполнить расчет и уточнение сегментации.

8       Сформировать и решить матричное уравнение.

9    *Иначе Если i = 2*

10       Выполнить расчет и уточнение сегментации.

11       Сформировать и решить матричное уравнение.

12       Сравнить матрицы СЛАУ, полученные на шагах 8 и 11 для нахождения изменившихся элементов.

13    *Иначе*

14       Выполнить расчет и уточнение сегментации.

15       Переформировать матрицу СЛАУ и решить матричное уравнение.

Поясним наиболее значимые шаги псевдокода. На первом шаге необходимо выбрать изменяемые параметры и задать диапазон их значений. Здесь же задается процент учащаемых сегментов $k$ и остальные параметры ЛП: геометрические параметры, значения относительных диэлектрических проницаемостей диэлектрических границ, число проводников. На втором шаге задается начальная сегментация (неравномерная сегментация с неизменным числом сегментов для каждой границы). Для достижения примерно равномерного разбиения проводников по осям координат выбирается начальная сегментация, при которой на всех границах исследуемой МПЛП, перпендикулярных оси абсцисс, отложено 3 сегмента, а перпендикулярных оси ординат – 40 сегментов. На третьем шаге производится предварительное решение матричного уравнения при среднем значении изменяемого параметра (медианном значении, если все элементы диапазона эквидистантны, среднем арифметическом значении в ином случае). На четвертом шаге в соответствии с методом I определяются сегменты, подлежащие разбиению. После этого в цикле производится расчет матрицы СЛАУ для каждого значения изменяемого параметра из диапазона его значений. На каждой итерации цикла сначала производится учащение исходной сегментации (шаги 7, 10 и 14) за счет деления сегментов (шаг 4) на два. При этом на первых двух итерациях происходит расчет всех элементов матрицы СЛАУ (шаги 8 и 11). На второй итерации дополнительно сравниваются элементы полученной матрицы и матрицы, полученной на шаге 8, для нахождения изменившихся элементов (шаг 12). На последующих итерациях за основу берется матрица СЛАУ, сформированная на шаге 8, с последующим пересчетом только её изменяющихся элементов, найденных на шаге 12.

Для иллюстрации принципа учета частичных изменений в матрице СЛАУ на рисунке 4 показана разница в портретах матрицы СЛАУ на первой и второй итерациях при изменении параметра $t$, где черным цветом показаны изменяемые элементы, а белым – неизменяемые. Здесь неизменяемые элементы составили 54% от их общего числа. Это свидетельствует о том, что затраты времени на непосредственный пересчет элементов матрицы СЛАУ могут быть сокращены более чем вдвое.

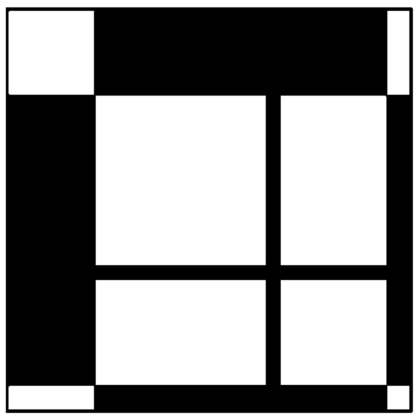

Рисунок 4 – Разница в портретах матрицы СЛАУ

на первой и второй итерациях при изменении $t$

Для демонстрации эффективности метода I использовано сравнение со случаем полного пересчета матрицы СЛАУ на каждой итерации (далее – метод II) на примере многовариантного анализа $m$-проводной ЛП при последовательном изменении её параметров $t$, $w$ и $\varepsilon_{r2}$. При вычислениях общие геометрические размеры МПЛП не изменялись. Параметры $t$, $w$ и $\varepsilon_{r2}$ изменялись в диапазоне ±14% (с шагом в 2%, всего 15 значений).

Результаты, полученные при $m = 6$, 8, 10 и 12, сведены в таблицы 9 – 11. В таблицах также приведено общее время расчета $t_{\text{tot}}$, требуемое на выполнение всех 15 итераций, среднее время $t_{\text{mid}}$ одной итерации и экономия времени относительно метода II.

Таблица 9 – Сравнение затрат времени методов I и II при изменении параметра $t$ для МПЛП 1 (рисунок 2)

| $m$ | Метод I | | Метод II | | $(t_{\text{tot}}^{\text{II}} - t_{\text{tot}}^{\text{I}})/t_{\text{tot}}^{\text{II}}$, % |
|---|---|---|---|---|---|
|  | $t_{\text{mid}}$, с | $t_{\text{tot}}$, с | $t_{\text{mid}}$, с | $t_{\text{tot}}$, с |  |
| 6 | 0,48 | 7,59 | 0,68 | 10,15 | 25 |
| 8 | 0,77 | 12,20 | 1,07 | 16,03 | 24 |
| 10 | 1,16 | 18,22 | 1,53 | 22,98 | 21 |
| 12 | 1,54 | 23,89 | 2,02 | 30,27 | 21 |

Таблица 10 – Сравнение затрат времени методов I и II при изменении параметра $w$ для МПЛП 1 (рисунок 2)

| $m$ | Метод I | | Метод II | | $(t_{\text{tot}}^{\text{II}} - t_{\text{tot}}^{\text{I}})/t_{\text{tot}}^{\text{II}}$, % |
|---|---|---|---|---|---|
|  | $t_{\text{mid}}$, с | $t_{\text{tot}}$, с | $t_{\text{mid}}$, с | $t_{\text{tot}}$, с |  |
| 6 | 0,60 | 9,06 | 0,68 | 10,15 | 11 |
| 8 | 0,89 | 13,64 | 1,07 | 16,03 | 15 |
| 10 | 1,33 | 20,11 | 1,53 | 22,98 | 12 |
| 12 | 1,96 | 29,12 | 2,02 | 30,27 | 4 |

Таблица 11 – Сравнение затрат времени методов I и II при изменении параметра $\varepsilon_{r2}$ для МПЛП 1 (рисунок 2)

| $m$ | Метод I | | Метод II | | $(t_{\text{tot}}^{\text{II}} - t_{\text{tot}}^{\text{I}})/t_{\text{tot}}^{\text{II}}$, % |
|---|---|---|---|---|---|
|  | $t_{\text{mid}}$, с | $t_{\text{tot}}$, с | $t_{\text{mid}}$, с | $t_{\text{tot}}$, с |  |
| 6 | 0,33 | 5,65 | 0,68 | 10,15 | 44 |
| 8 | 0,47 | 8,16 | 1,07 | 16,03 | 49 |
| 10 | 0,69 | 11,82 | 1,53 | 22,98 | 49 |
| 12 | 1,04 | 17,38 | 2,02 | 30,27 | 43 |

Как видно из таблиц 9 – 11, экономия времени расчета составила 4 – 49%. При этом наибольшая экономия достигается при изменении $\varepsilon_{r2}$. Стоит отметить, что, хотя метод I и требует дополнительного предварительного вычисления, при большом числе значений изменяемого параметра экономия времени существенно превышает накладные расходы на дополнительную итерацию, причем, чем шире диапазон изменяемого параметра, тем больше будет экономия времени. При одновременном изменении сразу нескольких параметров, требуется всего одно предварительное вычисление. Так, аналогичные результаты при одновременном изменении параметров $t$ и $\varepsilon_{r2}$ в тех же диапазонах приведены в таблице 12. Видно, что и в этом случае экономия времени достигает 23%.

Таблица 12 – Сравнение затрат времени методов I и II при изменении параметров $t$ и $\varepsilon_{r2}$ для МПЛП 1 (рисунок 2)

| $m$ | Метод I | | Метод II | | $(t_{tot}^{II} - t_{tot}^{I})/t_{tot}^{II}$, % |
|---|---|---|---|---|---|
| | $t_{mid}$, с | $t_{tot}$, с | $t_{mid}$, с | $t_{tot}$, с | |
| 6 | 0,53 | 8,30 | 0,68 | 10,15 | 18 |
| 8 | 0,76 | 12,29 | 1,07 | 16,03 | 23 |
| 10 | 1,15 | 18,22 | 1,53 | 22,98 | 21 |
| 12 | 1,60 | 24,66 | 2,02 | 30,27 | 19 |

Многовариантный анализ также апробирован на 10-проводной ЛП, поперечное сечение которой приведено на рисунке 3. Параметры $t$ и $\varepsilon_r$ изменялись в тех же диапазонах как по отдельности, так и совместно. Полученные результаты приведены в таблице 13.

Таблица 13 – Сравнение затрат времени методов I и II для МПЛП 2 (рисунок 3)

| Изменяемый параметр | Метод I | | Метод II | | $(t_{tot}^{II} - t_{tot}^{I})/t_{tot}^{II}$, % |
|---|---|---|---|---|---|
| | $t_{mid}$, с | $t_{tot}$, с | $t_{mid}$, с | $t_{tot}$, с | |
| $t$ | 0,83 | 12,91 | 1,14 | 17,16 | 25 |
| $\varepsilon_r$ | 0,50 | 8,68 | | | 49 |
| $t$ и $\varepsilon_r$ | 0,82 | 12,87 | | | 25 |

**Заключение**

Предложен метод неравномерной сегментации микрополосковых МПЛП с торцевой связью, показавший физичные результаты при вычислении емкостной матрицы. При этом его использование при многовариантном анализе МПЛП позволяет достичь экономии времени вычисления до 49% за счет учета частичных изменений в матрице СЛАУ, что подтверждает его эффективность. В дальнейшем целесообразна его апробация при решении оптимизационных задач МПЛП.

**Литература**

**Приложение А**

Таблица А.1 – Значения (пФ/м) первой строки матрицы **C** при $t = 0{,}018$ мм и изменении $n$, полученные в GNU Octave для МПЛП 1 (рисунок 2)

| $n$ | $C_{11}$ | $|C_{12}|$ | $|C_{13}|$ | $|C_{14}|$ | $|C_{15}|$ | $|C_{16}|$ | $|C_{17}|$ | $|C_{18}|$ |
|---|---|---|---|---|---|---|---|---|
| 1 | 101,00 | 15,80 | 0,60 | 7,97e–2 | 1,22e–2 | 2,05e–3 | 4,55e–4 | 2,46e–4 |
| 3 | 101,00 | 16,70 | 0,64 | 8,52e–2 | 1,29e–2 | 2,00e–3 | 3,22e–4 | 1,20e–4 |
| 5 | 101,00 | 16,90 | 0,64 | 8,60e–2 | 1,31e–2 | 2,02e–3 | 3,21e–4 | 1,28e–4 |
| 7 | 101,00 | 16,90 | 0,64 | 8,63e–2 | 1,31e–2 | 2,03e–3 | 3,23e–4 | 1,29e–4 |
| 9 | 101,00 | 16,90 | 0,65 | 8,65e–2 | 1,31e–2 | 2,04e–3 | 3,24e–4 | 1,29e–4 |

Таблица А.2 – Значения (пФ/м) первой строки матрицы **C** при $t = 0{,}05$ мм и изменении $n$, полученные в GNU Octave для МПЛП 1 (рисунок 2)

| $n$ | $C_{11}$ | $|C_{12}|$ | $|C_{13}|$ | $|C_{14}|$ | $|C_{15}|$ | $|C_{16}|$ | $|C_{17}|$ | $|C_{18}|$ |
|---|---|---|---|---|---|---|---|---|
| 1 | 120,00 | 24,40 | 0,53 | 5,55e–2 | 1,01e–2 | 3,76e–3 | 2,31e–3 | 1,86e–3 |
| 3 | 121,00 | 26,90 | 0,49 | 4,45e–2 | 4,94e–3 | 7,65e–4 | 2,62e–4 | 3,19e–4 |
| 5 | 121,00 | 27,30 | 0,49 | 4,51e–2 | 4,79e–3 | 5,68e–4 | 1,07e–4 | 1,62e–4 |
| 7 | 121,00 | 27,50 | 0,50 | 4,53e–2 | 4,77e–3 | 5,26e–4 | 7,40e–5 | 1,45e–4 |
| 9 | 121,00 | 27,50 | 0,50 | 4,54e–2 | 4,77e–3 | 5,14e–4 | 6,46e–5 | 1,45e–4 |

Таблица А.3 – Значения (пФ/м) первой строки матрицы **C** при $t = 0{,}07$ мм и изменении $n$, полученные в GNU Octave для МПЛП 1 (рисунок 2)

| $n$ | $C_{11}$ | $|C_{12}|$ | $|C_{13}|$ | $|C_{14}|$ | $|C_{15}|$ | $|C_{16}|$ | $|C_{17}|$ | $|C_{18}|$ |
|---|---|---|---|---|---|---|---|---|
| 1 | 136,00 | 32,50 | 0,67 | 5,23e–2 | 1,33e–2 | 7,06e–3 | 4,83e–3 | 4,15e–3 |
| 3 | 134,00 | 32,90 | 0,68 | 5,69e–1 | 4,90e–1 | 1,66e–2 | 4,59e–1 | 5,12e–1 |
| 5 | 133,00 | 33,90 | 0,37 | 2,46e–2 | 2,15e–3 | 3,35e–4 | 1,47e–4 | 2,86e–4 |
| 7 | 132,00 | 33,90 | 0,37 | 2,46e–2 | 2,01e–3 | 2,22e–4 | 6,46e–5 | 2,08e–4 |
| 9 | 132,00 | 34,00 | 0,37 | 2,47e–2 | 1,97e–3 | 1,83e–4 | 3,65e–5 | 1,92e–4 |

Таблица А.4 – Значения (пФ/м) первой строки матрицы **C** при $t = 0{,}105$ мм и изменении $n$, полученные в GNU Octave для МПЛП 1 (рисунок 2)

| $n$ | $C_{11}$ | $|C_{12}|$ | $|C_{13}|$ | $|C_{14}|$ | $|C_{15}|$ | $|C_{16}|$ | $|C_{17}|$ | $|C_{18}|$ |
|---|---|---|---|---|---|---|---|---|
| 1 | 153,00 | 56,80 | 3,54 | 3,74e–1 | 8,90e–2 | 4,76e–2 | 3,44e–2 | 4,04e–2 |
| 3 | 149,00 | 44,80 | 0,27 | 1,50e–2 | 4,03e–3 | 2,35e–3 | 1,65e–3 | 1,98e–3 |
| 5 | 153,00 | 44,90 | 0,19 | 8,16e–3 | 2,21e–3 | 1,40e–3 | 1,07e–3 | 2,33e–3 |
| 7 | 154,00 | 45,10 | 0,19 | 6,27e–3 | 4,10e–4 | 1,35e–4 | 1,03e–4 | 6,02e–4 |
| 9 | 154,00 | 45,20 | 0,19 | 6,24e–3 | 2,79e–4 | 3,41e–5 | 2,32e–5 | 4,01e–4 |